\documentclass[11pt]{article}  

\usepackage{epsfig}
\usepackage{amsfonts, amsmath, amssymb, color, mathrsfs, mathabx}  
 \usepackage{amsthm}
 \usepackage{placeins}
 
\usepackage{graphicx}
\usepackage{dcolumn}
\usepackage{bm}
\usepackage{caption}
\usepackage[utf8]{inputenc}
\usepackage[T1]{fontenc}
\usepackage{mathptmx}
\usepackage{etoolbox}
\usepackage{algorithm}
\usepackage{algpseudocode}
\usepackage{blindtext}

\usepackage{import}
\usepackage{subfiles} 
\usepackage{xcolor}
\DeclareCaptionFormat{cont}{#1 (cont.)#2#3\par}

\begin{document}

\title{A detailed analysis of the origin of deep-decoupling oscillations}
\author{John Bailie, Henk A. Dijkstra and Bernd Krauskopf}
\date{May 2024}
\maketitle

\begin{abstract}
\noindent
The variability of the strength of the Atlantic Meridional Overturning Circulation is influenced substantially by the formation of deep water in the North Atlantic. In many ocean models, so-called deep-decoupling oscillations have been found, whose timescale depends on the characteristics of convective vertical mixing processes.  Their precise origin and sensitivity to the representation of mixing have remained unclear so far. To study this problem, we revisit a conceptual Welander model for the evolution of temperature and salinity in two vertically stacked boxes for surface and deep water, which interact through diffusion and/or convective adjustment. The model is known to exhibit several types of deep-decoupling oscillations, with phases of weak diffusive mixing interspersed with strong convective mixing, when the switching between them is assumed to be instantaneous. We present a comprehensive study of oscillations in Welander's model with non-instantaneous switching between mixing phases, as described by a smooth switching function. A dynamical systems approach allows us to distinguish four types of oscillations, in terms of their phases of diffusive versus convective mixing, and to identify the regions in the relevant parameter plane where they exist. The characteristic deep-decoupling oscillations still exist for non-instantaneous switching, but they require switching that is considerably faster than needed for sustaining oscillatory behaviour. Furthermore, we demonstrate how a gradual freshwater influx can lead to transitions between different vertical mixing oscillations. Notably, the convective mixing phase becomes shorter and even disappears, resulting in long periods of much reduced deep water formation. The results are relevant for the interpretation of ocean-climate variability in models and (proxy) observations.
\end{abstract}

\section{Introduction}

Decades of studies of the North Atlantic Ocean circulation with Ocean General Circulation Models \cite[]{Mikolajewicz1990, Weaver1993, Delworth2000, LeBars2016, Lohmann2024} have shown that convective mixing, in particular, in the Greenland Sea and the Labrador Sea, is involved in the variability of the Atlantic Meridional Overturning Circulation (AMOC) on timescales of decades to millennia. This type of strong mixing between the top layer and deeper parts of the ocean modifies density gradients locally in the respective geographic regions, which lead to a large scale AMOC response via associated pressure gradients. 

A particular example are the so-called Dansgaard-Oeschger (D-O) events, which display variability \cite[]{menviel2020ice} on a millennial timescale. Their modelling has shown that the variability arises through the interaction of ocean-atmosphere, sea-ice, and convective mixing, which couples changes in atmospheric buoyancy fluxes to the AMOC variations. It has been speculated that so-called deep-decoupling oscillations may describe the D-O events.  These oscillations are characterized by a short phase of convective mixing in the North Atlantic followed by a very long deep-decoupling phase. During the deep-decoupling phase, mixing is diffusive and the AMOC re-adjusts, leading to overall low-frequency variations \cite[]{winton1993thermohaline}. For this variability, the AMOC's response to surface buoyancy flux anomalies is sensitive to the representation of the vertical mixing component \cite[]{Weijer2003, Vettoretti2022}. 

Conceptual models have played an important role in investigating the sensitivity of the AMOC to representations of vertical mixing \cite[]{PCessi1996, ColindeVerdiere2007}. The model by \cite{welander1986thermohaline} for vertical mixing in the North Atlantic takes a central place here, because it exhibits periodic solutions with the relevant deep-decoupling properties \cite[]{winton1995energetics, winton1993thermohaline, sima2004younger}. More precisely, it takes the form of a differential equation for the evolution of temperature and salinity, and its oscillations are characterized by two phases: a deep-coupling phase, where convective mixing between the surface and deep ocean in the North Atlantic strongly `couples' the surface and deep ocean, and a deep-decoupling phase when a polar halocline prevents convective mixing. The \cite{welander1986thermohaline} model has been analysed extensively by \cite{PCessi1996} for the case of instantaneous switching between convective and diffusive mixing phases, as represented by a step function. \cite{leifeld2018non} studied the mathematical mechanism of how the periodic solution (dis)appears in this piecewise-smooth model. A comprehensive bifurcation analyses of all dynamical regimes of this limiting case of the Welander model has been presented in \cite{bailie2024bifurcation}.

While instantaneous switching between different phases is an interesting limiting case, most models will feature a continuous and smooth switching function to represent switching behaviour. In the Welander model, which is introduced in detail in Sec.~\ref{sec:introduceModel}, switching is realised by a $\tanh$ function with a switching timescale parameter $\varepsilon$. This switching function becomes a step function for $\varepsilon = 0$, but otherwise the model is smooth. The existence of periodic solutions for $\varepsilon > 0$ was also studied in \cite{bailie2024bifurcation}; the arguably surprising result is that the Welander model supports oscillations only when the switching between convective and diffucsive mixing phases is sufficiently fast.

We present here an investigation into the nature of the oscillations of the Welander model for positive $\varepsilon$ in terms of the existence of deep-coupling phases where mixing is convective, and deep-decoupling phases characterised by weaker diffusive mixing. We take a mathematical approach grounded in dynamical system theory for systems of smooth differential equations \cite[]{kuznetsov1998elements, strogatz2000} in combination with the continuation software {\sc AUTO} \cite[]{doedel2007auto}. We define the boundaries of the switching zone for any $\varepsilon \geq 0$, and this enables us to determine when the periodic solution has a component in the zones of deep coupling and/or of deep decoupling. For the limiting case $\varepsilon = 0$, we find  only oscillations with both a deep-coupling and a deep-decoupling phase \cite[]{PCessi1996,leifeld2018non,bailie2024bifurcation}, which we refer to as \emph{Welander oscillations}. For $\varepsilon > 0$, however, these are no longer the only type of oscillation. We present a bifurcation analysis that distinguishes different sub-regions in the relevant parameter plane (of virtual salinity flux and density threshold), where the corresponding oscillation has qualitatively different features. In particular, we find that the existence of Welander deep-(de)coupling oscillations requires even faster switching than needed to sustain periodic solutions. 

The paper is organised as follows. Section~\ref{sec:introduceModel} introduces Welander's model, and this is followed in Sect.~\ref{sec:non_smooth} by a brief presentation of the oscillatory regime for the piecewise-smooth limiting case $\varepsilon = 0$. In Sec.~\ref{sec:zones}, we define the boundary of the switching zone and the zones of deep-(de)coupling for $\varepsilon > 0$. Section~\ref{sec:BifDiag} then discusses the different types of oscillations and where they can be found for the representative value $\varepsilon = 0.009$.  We start by showing that Welander deep-(de)coupling oscillations still exist, and show in Sec.~\ref{sec:BifDiag}\ref{sec:tangencies} that they are lost when the underlying periodic orbit becomes tangent to the boundary of the switching zone. In Sec.~\ref{sec:BifDiag}\ref{sec:existence}, we then present bifurcation diagrams that show for which combinations of virtual salinity flux and density threshold Welander and other types of oscillations exist; their properties in terms of different oscillation phases are detailed in Sec.~\ref{sec:BifDiag}\ref{sec:phaseportraits}. By way of an interpretation of these results, Sec.~\ref{sec:drift} discusses how the nature of the observed oscillation changes gradually when the virtual salinity flux slowly decreases, as would be the case with increasing influx of meltwater from the Greenland ice sheet. Section~\ref{sec:disapear_occ} then shows how the bifurcation diagram in the parameter plane changes with increasing $\varepsilon$, as the region of periodic solutions shrinks and then disappears. The final Sec.~\ref{sec:conclusions} discusses our findings and points out directions for future investigations.

\section{Conceptual mixing model}
\label{sec:introduceModel}

\begin{figure}[t!]
	\centering
	\includegraphics{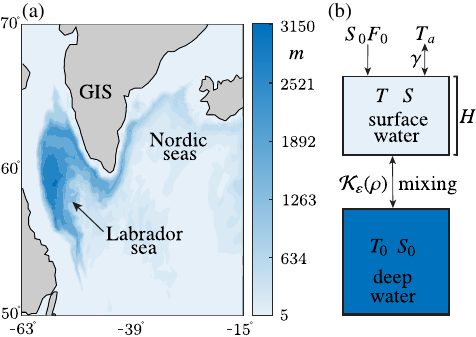}
	\caption{Deep water formation in the North Atlantic and its modelling. Panel (a) presents the Greenland Ice Sheets (GIS), and the Labrador and Nordic seas with the average mixing depth in meters during January 1958 using ORAS5 global ocean reanalysis data \cite[]{Copernicus}; here, dark blue regions illustrate deeper mixing depths, indicating areas of active deep water formation. Panel (b) is a schematic of the two-box Welander model as described by Eqs.~\eqref{eq:dimWelander}.}
	\label{fig:schematic}
\end{figure}

The \cite{welander1986thermohaline} model comprises of two vertically stacked boxes representing the surface and deep ocean near deep water formation sites in the subpolar North Atlantic, as shown in Fig.~\ref{fig:schematic}(a). It describes the evolution of salinity $S$ and temperature $T$ in the surface box of height $H$ that is coupled to a deep water box with prescribed and fixed temperature $T_0$ and salinity $S_0$. Figure \ref{fig:schematic}(b) presents a schematic of this model. Temperature and salinity between the surface and bottom boxes interact according to a mixing process with the exchange function $K_\varepsilon(\rho)$ that depends on the density $\rho$ at the surface relative to that of the bottom box. Moreover, the salinity $S$ in the surface box experiences the virtual salinity flux $S_0F_0$, while the temperature $T$ is relaxing to the atmospheric temperature $T_a$ at rate $\gamma$. The overall model takes the form 
\begin{equation}
\begin{aligned} 
	&\frac{dT}{dt} = -\gamma(T - T_a) - (T - T_0)K_\varepsilon(\rho),
	\\
	&\frac{dS}{dt} = \frac{F_0S_0}{H} - (S - S_0)K_\varepsilon(\rho).
	\label{eq:dimWelander}
\end{aligned}
\end{equation}
Here, we parameterize the vertical mixing process between the two boxes with enhanced diffusion, as also done in earlier studies \cite[]{vellinga1998multiple, den2011spurious}; specifically, the convective exchange function
\begin{align}
	K_\varepsilon(\rho)  =  k_1 + \frac{1}{2}(k_2 - k_1) \left(1 + \tanh{\left(\frac{\rho - \rho_0 - g^*}{\varepsilon}\right)}\right),
	\label{eq:conv_ex}
\end{align}
smoothly transitions between diffusive and convective vertical mixing phases, given by vertical mixing coefficients $k_1$ and $k_2$, where $0<k_1<k_2$. Further, $g^{*}$ is the density threshold that controls the transition between the two types of mixing. When the density difference $\rho-\rho_0$ between the surface and bottom box exceeds $g^{*}$, the water column becomes unstably stratified. Consequently, denser surface water overlies a lighter deep water, promoting strong, vigorous convective mixing between the boxes, represented by $k_2$. In contrast, when $\rho-\rho_0$ falls below $g^{*}$, the water column is stably stratified, resulting in weak, diffusive mixing represented by $k_1$. Note that the slope of $K_\varepsilon(\rho)$ at $\rho - \rho_0 = g^*$ is $1/\varepsilon$. The smaller $\varepsilon$, the steeper is this slope and, hence, the faster the transition, which is why we refer to $\varepsilon$ as the switching timescale parameter. Indeed, for $\varepsilon = 0$ the transition is `infinitely fast', that is, instantaneous.

To complete the model and express $K_\varepsilon(\rho)$ as a function of $S$ and $T$, we follow standard practice in many conceptual box models \cite[]{welander1986thermohaline, PCessi1996, weijer2019stability} and employ a linear equation of state 
\begin{align}
	\frac{\rho}{\rho_0} = 1 + \alpha_S(S - S_0) - \alpha_T(T - T_0).
	\label{eq:def_rho}
\end{align}
Here, $\alpha_S$ and $\alpha_T$ are the saline expansion and thermal compression coefficients, respectively, and the bottom box is taken as a prescribed reference state to measure the relative dynamics of the surface box.

The analysis of system~\eqref{eq:dimWelander} with $K_\varepsilon(\rho) = K_\varepsilon(S,T)$ as given by~\eqref{eq:conv_ex} and~\eqref{eq:def_rho} is greatly simplified by non-dimensionalising  \cite[]{PCessi1996, bailie2024bifurcation}, which reduces the number of parameters. The resulting model is given by 
\begin{equation}
	\begin{aligned}
		&\dot x = 1 - x - \mathcal{K}_\varepsilon(x,y)x,
		\\
		&\dot y = \mu - \mathcal{K}_\varepsilon(x,y)y,
	\end{aligned}
	\label{eq:nondimWel}
\end{equation}
for the non-dimensional temperature $x$ and salinity $y$; here, the dot represents the derivative with respect to the rescaled time and $\mu$ is the rescaled virtual salinity flux. Additionally, the convective exchange function now takes the form 
\begin{equation}
	\begin{aligned}
		&\mathcal{K}_\varepsilon(x,y) = \kappa_1 + \frac{1}{2}(\kappa_2 - \kappa_1) \left(1 + \tanh{\left(\frac{y - x - \eta}{\varepsilon}\right)}\right),
	\end{aligned}
	\label{eq:convExchang}
\end{equation}
where $\kappa_1$ and $\kappa_2$ are the non-dimensional vertical mixing coefficients and $\eta$ is the non-dimensional density threshold. Conveniently, the non-dimensional density takes the simple form $\rho = y - x$, so that we have $\mathcal{K}_\varepsilon(\rho) = \mathcal{K}_\varepsilon(x,y)$. 

We refer to system~\eqref{eq:nondimWel} with~\eqref{eq:convExchang} as the adjusted Welander model, and it is our main object of study. Specifically, we fix the vertical mixing coefficients at $\kappa_1 = 0.1$ and $\kappa_2 = 1.0$, as in \cite[]{bailie2024bifurcation}, and investigate the location in the $(\mu,\eta)$-plane of oscillations with both deep-decoupling and deep-coupling phases, and how their existence depends on the switching timescale parameter $\varepsilon$. To this end, we make use of the fact that the adjusted Welander model is an example of a `switched' slow-fast system \cite[]{wechselberger2020geometric, kristiansen2021relaxation}. More specifically, the convective exchange function $\mathcal{K}_\varepsilon(\rho)$ in~\eqref{eq:convExchang} is a switch from $\kappa_1$ for small $\rho$ to $\kappa_2$ for large $\rho$; the switching happens over a $\rho$-interval around $\rho = \eta$, whose width shrinks to $0$ with $\varepsilon$.

\section{Welander oscillations of the piecewise-smooth limit}
\label{sec:non_smooth}

\begin{figure*}[t!]
    \centering
    \includegraphics{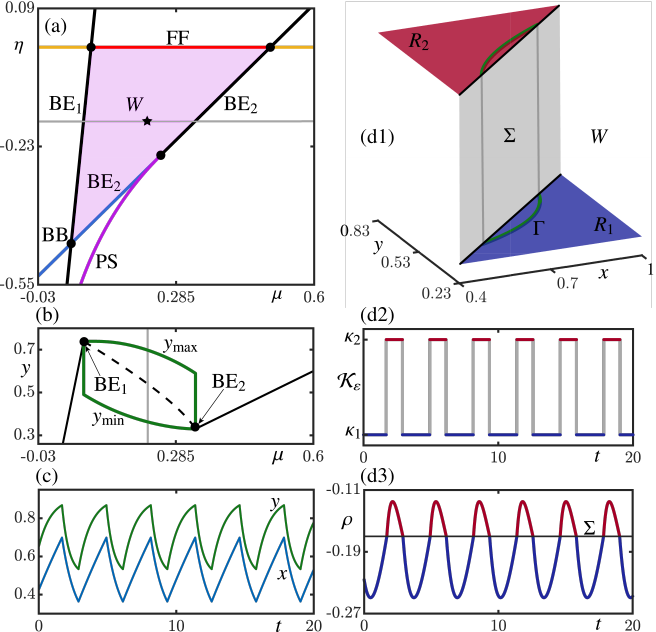}
    \caption{\label{fig:BDnonsmooth} Bifurcation diagrams and representations of the stable periodic solution of system~\eqref{eq:nondimWel} with $\varepsilon=0$. Panel~(a) is the bifurcation diagram in the $(\mu,\eta)$-plane with curves $\mathrm{BE_1}$, $\mathrm{BE_2}$, $\mathrm{BB}$, and $\mathrm{PS}$ of discontinuity-induced bifurcations (coloured as in \cite[]{bailie2024bifurcation}) that bound the region $W$ of Welander oscillations (shaded purple). Panel~(b) shows the one-parameter bifurcation diagram in $\mu$ for fixed $\eta= -0.17$, with branches of equilibria (black; solid when stable, dashed when unstable) and the branch of periodic solutions represented by their minimum and maximum values $y_{\rm min}$ and $y_{\rm max}$ (green curves). Panel~(c) shows the time series of temperature $x$ (blue) and salinity $y$ (green) for $(\mu,\eta) = (0.219, -0.17)$, at the marker in panel (a1). Panel~(d1) shows the corresponding stable periodic solution $\Gamma$ (green) on the piecewise constant graph of $\mathcal{K}_\varepsilon$ over the $(x, y)$-plane; here, the graph is coloured blue where $\mathcal{K}_\varepsilon = \kappa_1$ over $R_1$ and red where $\mathcal{K}_\varepsilon = \kappa_2$ over $R_2$, which are connected by the instantaneous vertical switching plane $\Sigma$ (grey). Panels~(d2) and (d3) are the correspondingly coloured time series of $\mathcal{K}_\varepsilon$ and $\rho$ along $\Gamma$, respectively.}
\end{figure*}

In the limit $\varepsilon = 0$, the switching between $\kappa_1$ and $\kappa_2$ is instantaneous and system~\eqref{eq:nondimWel} is no longer smooth; rather, it is now a piecewise linear system with a discontinuity along the line where $\rho = y - x = \eta$. As the starting point of our investigation, we recall and present where in the $(\mu,\eta)$-plane one finds a piecewise smooth periodic orbits, corresponding to Welander oscillations (see  \cite{bailie2024bifurcation} for a complete bifurcation analysis of this limiting case). 

For $\varepsilon = 0$, system~\eqref{eq:nondimWel} with~\eqref{eq:convExchang} is known as a planar Filipov system \cite[]{bernardo2008piecewise, filippov2013differential, guardia2011generic} with the two adjacent regions or zones
\begin{align*}
    R_1 &= \{(x,y)\in\mathbb{R}^2, \ y < x + \eta\},
    \\ 
    R_2 &= \{(x,y)\in\mathbb{R}^2, \ y > x + \eta\},
\end{align*}
which are separeted by the switching line
\begin{align*}
    \Sigma = \{(x,y)\in\mathbb{R}^2, \ y = x + \eta\}. 
\end{align*}
Zone $R_1$ represents a stably stratified water column with diffusive mixing and,  conversely, zone $R_2$ represents an unstably stratified water column with associated vigorous convective mixing that couples the surface and deep water boxes. The dynamics in $R_1$ is determined by setting $\mathcal{K}_\varepsilon(\rho) \equiv \kappa_1$ in~\eqref{eq:nondimWel}, and that in $R_2$ by setting $\mathcal{K}_\varepsilon(\rho) \equiv \kappa_2$. Both of these systems are linear, and this implies that any periodic solutions of the planar Filipov system~\eqref{eq:nondimWel} with $\varepsilon = 0$ must consist of a segment in zone $R_1$ and a segment in zone $R_2$ and, thus, cross the switching line $\Sigma$ twice. 

Figure~\ref{fig:BDnonsmooth} illustrates that this type of oscillating solution exists, namely for parameter combinations of $\mu$ and $\eta$ in the shaded parameter region labeled W of the $(\mu,c)$-plane shown in panel~(a). Along the curves $\mathrm{FF}$, $\mathrm{BE}_1$, $\mathrm{BB}$ and $\mathrm{BE_2}$ that bound this region one finds discontinuity-induced bifurcations \cite[]{filippov2013differential, di2008discontinuity, bernardo2008piecewise} that lead to the emergence or disappearance of the periodic solution \cite[]{bailie2024bifurcation}. When $\eta$ is fixed at a correpsonding value, a stable periodic solution $\Gamma$ exist for $\mu$-values between the discontinuity-induced bifurcations $\mathrm{BE_1}$ and $\mathrm{BE_2}$. This is illustrated in Fig.~\ref{fig:BDnonsmooth}(b) for $\eta = -0.17$, where $\Gamma$ is represented by its maximal and minimal $y$-values; note that system~\eqref{eq:nondimWel} with $\varepsilon = 0$ converges to stable equilibria outside the range of the stable periodic solution. Panel~(c) shows the stable oscillating solution for $(\mu,\eta) = (0.219, -0.17)$ as the time series of temperature $x$ and salinity $y$. Notice the piecewise smooth nature of these Welander oscillations, with corners (discontinuities of their derivative) at the maxima and minima; they are indeed consistent with  oscillations found in earlier studies \cite[]{PCessi1996,welander1986thermohaline}. 

To clarify which parts of the stable periodic orbit are in $R_1$ and in  $R_2$, respectively, we show $\Gamma$ in Fig.~\ref{fig:BDnonsmooth}(d1) on the piecewise-constant graph of the convective exchange function $\mathcal{K}_\varepsilon$ with $\varepsilon = 0$; the switching line $\Sigma$ is represented here by a vertical plane connecting the plateaux of constant $\kappa_1$ over $R_1$ and $\kappa_2$ over $R_2$. As is illustrated also by the time series of $\mathcal{K}_\varepsilon$ in Fig.~\ref{fig:BDnonsmooth}(d2), the periodic solution $\Gamma$ indeed features an alternation of deep-decoupling and deep-coupling phases of considerable lengths, with instantaneous switching at $\Sigma$ between them. This oscillating behaviour is represented in panel~(d3) by the time series of the density $\rho$ in the surface box. Notice the non-smooth nature of the oscillation whenever the threshold $\rho = \eta$ defining $\Sigma$ is crossed. The density is above this threshold during the deep-coupling phase, where the mixing is convective and vigorous. This phase ends abruptly when $\Gamma$ crosses $\Sigma$, which results in an immediate switch to a stable stratification of the water column. This is the beginning of the deep-decoupling phase with only weak, diffusive mixing, which ends when $\Gamma$ crosses $\Sigma$ again and there is, hence, an immediate switch back to the deep-coupling phase. We remark that the deep-decoupling phase of $\Gamma$ is longer when $\mu$ is closer to the boundary $\mathrm{BE}_1$ in Fig.~\ref{fig:BDnonsmooth}(a) and~(b). Conversely, higher values of $\mu$ near $\mathrm{BE_2}$ lead to a longer deep-coupling phase.

\section{Switching zone of the smooth model} 
\label{sec:zones}

Figure~\ref{fig:BDnonsmooth} shows that any periodic orbit of system~\eqref{eq:nondimWel} with $\varepsilon = 0$ necessarily has both a deep-coupling phase and a deep-decoupling phase, which is the characterising property of Welander oscillations. We now address the question whether and where Welander oscillations can be found when $\varepsilon > 0$. Our starting point is the analysis in \cite{bailie2024bifurcation}, which showed that there exists a region of the $(\mu,\eta)$-plane where a stable periodic orbit $\Gamma$ exists, provided $\varepsilon$ is sufficiently small, namely below $\varepsilon^* \approx 0.147$. However, no statement has been made previously as to when this periodic orbit is actually a Welander oscillation. 

For non-zero $\varepsilon > 0$ the transition between the  deep-decoupling phase near $\kappa_1$ and the deep-decoupling phase near $\kappa_2$ is no longer instantaneous. Instead, it occurs within an intermediate (narrow) switching zone of the $(x, y)$-plane that corresponds to the steep part of the now smooth convective exchange function $\mathcal{K}_\varepsilon$ from~\eqref{eq:convExchang}. To make this idea precise, one needs to decide when the switching begins and ends. This can be done in different ways, such as by setting threshold values for $\mathcal{K}_\varepsilon$ or by considering a percentage deviation from $\kappa_1$ and $\kappa_2$. 

We adopt here a geometric approach and define the switching part to lie in between the two points of maximal curvature of the graph of $\mathcal{K}_\varepsilon$ as a function of the density $\rho$. These points can be found from~\eqref{eq:convExchang} as the zeros of the derivative of the curvature along $\mathcal{K}_\varepsilon(\rho)$; their $\rho$-values are given as the roots $\rho^{\pm}$ of the function
\begin{equation}
\begin{aligned}
G(\rho) &= 8\varepsilon^2- 4(\kappa_1 - \kappa_2)
 \\  &+\left(9\varepsilon^2+8(\kappa_1-\kappa_2)^2\right)\cosh\left(\frac{2(\rho-\eta)}{\varepsilon}\right)
 \\
 &-\varepsilon^2\cosh\left(\frac{6(\rho-\eta)}{\varepsilon}\right),
	\label{eq:max_curve_root}
\end{aligned}
\end{equation}
and the corresponding $\kappa$-values are $L^{\pm} = \mathcal{K}_\varepsilon(\rho^{\pm})$ according to \eqref{eq:convExchang} with $\rho = y - x$.

\begin{figure}
    \centering
    \includegraphics{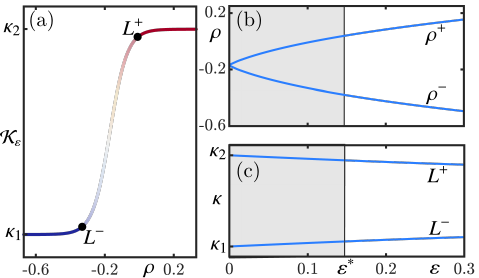}
    \caption{\label{fig:KappaRhovsEpsilon} The maximum curvature points $L^{\pm}$ of the convective exchange function $\mathcal{K}_\varepsilon(\rho)$. Panel~(a) shows them for $\varepsilon = 0.02$ and $\eta = -0.17$ on the graph of $\mathcal{K}_\varepsilon(\rho)$, which is coloured blue up to $L^{-}$ and red from $L^{+}$, with transitional, fainter colouring in between. Panels~(b) and~(c) show, respectively, the $\rho$-values $\rho^{\pm}$ and corresponding $\kappa$-values of $L^{\pm}$ as functions of $\varepsilon$ over the range $[0, 0.3]$.}
\end{figure}

Figure~\ref{fig:KappaRhovsEpsilon} illustrates the location of the maximum curvature points, which we also denote $L^{\pm}$ for notational convenience. Panel~(a) displays $\mathcal{K}_\varepsilon(\rho)$ for $\varepsilon = 0.02$ with $\eta = -0.17$, with $L^{-}$ and $L^{+}$ clearly lying where the curvature of this graph is largest. The color scheme indicates that $\mathcal{K}_\varepsilon(\rho)$ is near $\kappa_1$ to the left of $L^{-}$ and near $\kappa_2$ to the right of $L^{+}$, with a rapid but smooth switch between the two over the $\rho$-range in between $L^{-}$ and $L^{+}$. Panels~(b) and~(c) of Fig.~\ref{fig:KappaRhovsEpsilon} show how the points $L^{\pm}$ depend on $\varepsilon$ over a considerable range, which includes the value $\varepsilon^{*} \approx 0.147$ beyond which system~\eqref{eq:nondimWel} is no longer able to exhibit oscillations \cite[]{bailie2024bifurcation}. In the limit $\varepsilon = 0$ the associated $\rho$-values coincide at $\rho^{\pm} = \eta$ and the $\kappa$-values are $\kappa_1$ and $\kappa_2$, respectively. When $\varepsilon$ is increased from $0$, the distance in $\rho$ between $L^{-}$ and $L^{+}$ increases quite rapidly, while the corresponding $\kappa$-values remain very close to $\kappa_1$ and $\kappa_2$. 

Overall, Fig.~\ref{fig:KappaRhovsEpsilon} shows that it quite a natural choice, for any $\varepsilon \geq 0$, to consider the maximal curvature points $L^{\pm}$ as the boundary of the switching part of the convective exchange function $\mathcal{K}_\varepsilon(\rho)$. It results in the division of the $(x,y)$-plane of system~\eqref{eq:nondimWel} into the three distinct open regions
\begin{align}
	R_1 &= \{(x,y) \in \mathbb{R}^2, \  \mathcal{K}_\varepsilon(x,y) < L^{-} \}, \label{eq:regions_slowR1} \\
   S   &= \{(x,y,) \in \mathbb{R}^2, \ L^{-} < \mathcal{K}_\varepsilon(x,y) < L^{+}\},  \label{eq:regions_fastS} \\
	R_2 &= \{(x,y) \in \mathbb{R}^2, \  L^{+} < \mathcal{K}_\varepsilon(x,y)\},
 \label{eq:regions_slowR2} 
\end{align}
by the two lines along which $\mathcal{K}_\varepsilon(x,y) = L^{\pm}$. To avoid confusion with regions in the parameter $(\mu,\eta)$-plane we discuss in Sec.~\ref{sec:BifDiag}, we again rather speak of zones: the deep-decoupling zone $R_1$, the switching zone $S$, and the deep-coupling zone $R_2$.

For $\varepsilon = 0$ this division of the phase plane reduces to that into only the two zones $R_1$ and $R_2$, since then $S = \emptyset$ with $L^{\pm} = \Sigma$. For $\varepsilon > 0$, zone $R_1$ indeed still represents the deep-decoupling regime with $\mathcal{K}_\varepsilon(x,y)$ very near $\kappa_1$, and zone $R_2$ the deep-coupling regime with $\mathcal{K}_\varepsilon(x,y)$ very near $\kappa_2$. However, there is now a local timescale separation: in zones $R_1$ and $R_2$ the dynamics is much slower, compared to the fast switching between $\kappa_1$ and $\kappa_2$ in zone $S$. We remark that the dynamics in the intermediate region $S$ represents a slow-fast desingularization or ``blow up" with $\varepsilon$ of the switching line $\Sigma$ \cite[]{teixeira2012regularization, kuehn2015multiple}.

\section{Welander oscillations with smooth switching} 
\label{sec:BifDiag}

\begin{figure}[t!]
	\centering
	\includegraphics{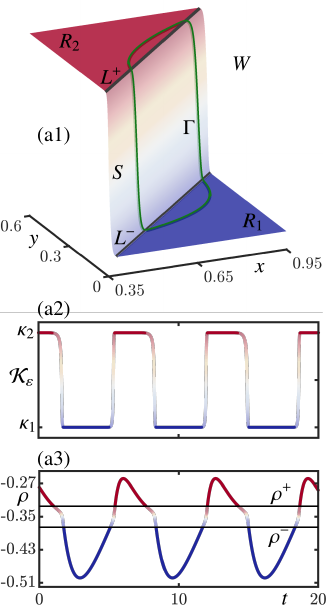}
	\caption{Welander oscillation at $(\mu,\eta) = (0.14, -0.3)$ of system~\eqref{eq:nondimWel} with $\varepsilon = 0.009$. Panel~(a1) shows the corresponding stable periodic solution $\Gamma$ (green) on the graph of $\mathcal{K}_\varepsilon$ (colored as in Fig.~\ref{fig:KappaRhovsEpsilon}(a)) with $L^{\pm}$, and  panels~(d2) and (d3) are the corresponding time series of $\mathcal{K}_\varepsilon$ and $\rho$. Compare with Fig.~\ref{fig:BDnonsmooth}(d1)--(d3).}
	\label{fig:Welander}
\end{figure} 

Welander oscillations can be found for sufficiently small $\varepsilon > 0$, and an example with $\varepsilon = 0.009$ is shown in Fig.~\ref{fig:Welander}. Panel~(a1) shows the underlying stable periodic orbit $\Gamma$ on the graph of $\mathcal{K}_\varepsilon$. It indeed features a deep-decoupling phase in $R_1$ and a deep-coupling phase in $R_2$, now with transitions through the switching zone $S$. The time series of $\mathcal{K}_\varepsilon$ in panels~(a2) illustrate that this oscillation indeed spends considerable amounts of time in $R_1$ and in $R_2$, respectively, with rapid transitions within $S$ on a faster timescale when $L^{\pm}$ are crossed. The density $\rho$ in panel~(a3) shows the same properties but is overall a more gradual oscillation: $\rho$ changes quite slowly during the deep-decoupled phase in $R_1$, where it has a minimum, and in the deep-coupling phase in $R_2$, where it has a maximum; the faster transtions through the switching zone $S$ occur in the range between $\rho^-$ and $\rho^+$. Comparison with Fig.~\ref{fig:BDnonsmooth}(d1)--(d3) shows that the Welander oscillation shown in  Fig.~\ref{fig:Welander} is the natural counterpart for $\varepsilon > 0$ of those of the piecewise smooth limit for $\varepsilon = 0$. Indeed, its characterizing feature is that the coexistence of these deep-(de)coupling phases generates extended periods of inactive and active deep water production; notice that the actual switching time is short here, namely about 8\% of the period of the oscillation.

\subsection{Tangencies of periodic solution with $L^{\pm}$}
\label{sec:tangencies}

\begin{figure*}[t!]
	\centering
    \includegraphics{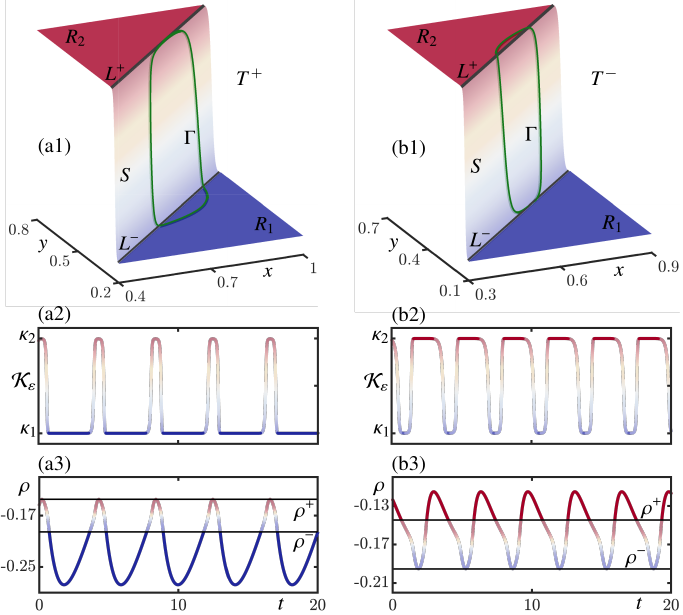}
	\caption{Tangencies of the periodic orbit $\Gamma$ of system~\eqref{eq:nondimWel} with $\varepsilon = 0.009$, shown as in Fig.~\ref{fig:Welander}. Panels~(a1)--(a3) show the tangency $T^{+}$ of $\Gamma$ with $L^{+}$ at $(\mu,\eta) = (0.1616, -0.17)$, and panels~(b1)--(b3) show its tangency $T^{-}$ with $L^{-}$ at $(\mu,\eta) = (0.3092, -0.17)$.}
	\label{fig:PPT}
\end{figure*}

The stable periodic orbit $\Gamma$ in Fig.~\ref{fig:Welander} can be followed or continued as the virtual salinity $\mu$ and the density threshold $\eta$ are varied. In the process, the lengths of its deep-decoupling phase in $R_1$ and its deep-coupling phase in $R_2$ change. In fact, either phase may shrink down to zero and disappear, and this happens when $\Gamma$ becomes tangent to $L^{-}$ or $L^{+}$. Figure~\ref{fig:PPT} illustrates two cases for fixed $\eta =-0.17$ in the style of Fig.~\ref{fig:Welander}. Panels~(a1)--(a3) of Fig.~\ref{fig:PPT} show the situation when $\Gamma$ is tangent to the upper threshold $L^{+}$, and we refer to this tangency as $T^{+}$. While the periodic orbit still crosses the entire switching zone $S$ to just reach $\rho^+$ at the maximum of $\rho$, it no longer enters the zone $R_2$ and, hence, $\Gamma$ no longer has a deep-coupling phase. Physically, the water column is almost sufficiently destabilized: changing the virtual salinity $\mu$ slightly will result in one of two changes: either a deep-coupling phase will form in zone $R_2$, evolving on a slower timescale and with a density sufficient for convective mixing to occur, or the solution will be completely contained within zones $S$ and $R_1$. Fig.~\ref{fig:PPT}(ab)--(b3) presents the converse scenario of the tangency $T^{-}$, where $\Gamma$ is tangent to the lower threshold $L^{-}$. Here the density has a minimum at $\rho^+$ and changing $\mu$ slightly will result in either the formation of a deep-decoupling phase, or the oscillation being restricted to zones $S$ and $R_2$.

\subsection{Existence region in the $(\mu,\eta)$-plane}
\label{sec:existence}

\begin{figure*}[t!]
	\centering
        \includegraphics{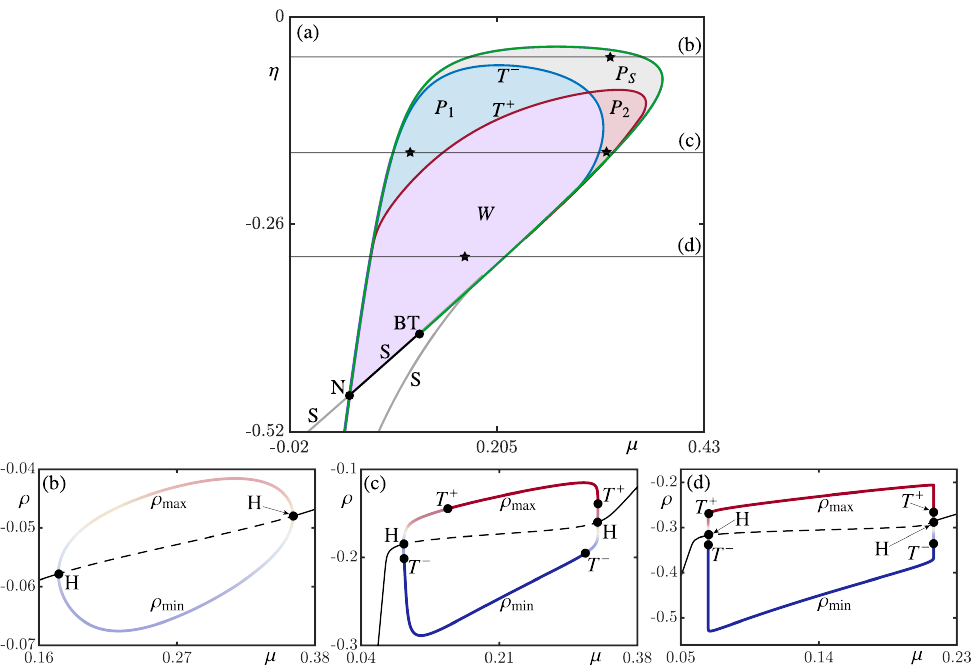}
	\caption{Bifurcation diagrams of system~\eqref{eq:nondimWel} with $\varepsilon = 0.009$. Panel~(a) shows the region in the $(\mu,\eta)$-plane where there exists a stable periodic orbit $\Gamma$; it is bounded by branches $\mathrm{H}$ of Hopf bifurcation (green) and $\mathrm{S}$ of saddle-node bifurcation (gray, and black when a SNIC), which meet at the shown points $\mathrm{BT}$ of Bogdanov-Takens and $\mathrm{N}$ of non-central saddle-node bifurcation. This region is divided by the tangency curves $T^{+}$ and $T^{-}$ into the four sub-regions $P_S$ (shaded gray), $P_1$ (shaded blue), $P_2$ (shaded red), and $W$ (shaded purple) of different types of oscillations. One-parameter bifurcation diagrams in $\mu$ along the accordingly labeled horizontal lines at $\eta = 0.05$, $\eta = -0.17$ and $\eta = -0.3$ are presented in panels~(b), (c) and~(d), respectively. Shown are the values of $\rho$ of branches of equilibria (black; solid when stable, dashed when unstable), and the minimum and maximum values $\rho_{\rm min}$ and $\rho_{\rm max}$ of $\Gamma$ between the two points $\mathrm{H}$, which are coloured by the value of $\mathcal{K}_\varepsilon(\rho)$ as in Fig.~\ref{fig:KappaRhovsEpsilon}(a); the tangency points $T^{+}$ and $T^{-}$ delimit the $\mu$-ranges of the deep-decoupling phase (blue) and the deep-coupling phase (red), respectively.}
	\label{fig:BD}
\end{figure*}

The stable periodic orbit $\Gamma$ of system~\eqref{eq:nondimWel} exists, for sufficiently small $\varepsilon$, in a region of the $(\mu,\eta)$-plane that is the continuation of the triangular region $W$ in Fig.~\ref{fig:BDnonsmooth}(a). Figure~\ref{fig:BD}(a) shows this region for $\varepsilon = 0.009$. It is bounded almost entirely by a curve $\mathrm{H}$ of Hopf bifurcation, all the way from the point $\mathrm{N}$ of non-central saddle-node bifurcation to the point $\mathrm{BT}$ of Bogdanov-Takens bifurcation (from where $\mathrm{H}$ emerges). Between $\mathrm{BT}$ and $\mathrm{N}$ the boundary is formed by the part of a curve $\mathrm{S}$ of saddle-node bifurcations, where the saddle-node lies on the periodic orbit $\Gamma$ (this is also known as a SNIC or SNIPER bifurcation \cite[]{kuznetsov1998elements,strogatz2000}). Additional bifurcation curves emerge from the points $\mathrm{N}$ and $\mathrm{BT}$, but they bound such narrow regions \cite[]{bailie2024bifurcation} that they are not of interest for the discussion here.

A crucial difference with the limiting case $\varepsilon = 0$ is that the stable periodic orbit $\Gamma$ for $\varepsilon > 0$ is not necessarily a Welander oscillation. In fact, when it is born at the Hopf bifurcation $\mathrm{H}$ the periodic orbit is initially very small and surrounds an equilibrium in the switching zone $S$ \cite[]{bailie2024bifurcation}; therefore, it is not a Welander oscillation near the curve $\mathrm{H}$. However, as $\Gamma$ grows, it has tangencies $T^{\pm}$ with $L^{\pm}$ and then enters the respective zones $R_1$ and/or $R_2$. In the $(\mu,\eta)$-plane these two types of tangencies from Fig.~\ref{fig:PPT} occur along the curves $T^{+}$ and $T^{-}$. These two curves can be found by continuing $\Gamma$ together with the condition that it is tangent to $L^{+}$ and $L^{-}$, respectively; the boundary value problem setup required for this computation is described briefly in Appendix~\ref{sec:computeT}. 

As Fig.~\ref{fig:BD}(a) shows, the curves $T^{\pm}$ divide the region of existence of $\Gamma$ into the four parameter sub-regions $P_S$, $P_1$, $P_2$ and $W$. As we discussed in Sec.~\ref{sec:tangencies}, when crossing $T^{+}$ and $T^{-}$, the periodic orbit $\Gamma$ gains or loses its deep-coupling or deep-decoupling phase. Hence, the bifurcation diagram in the $(\mu,\eta)$-plane provides a complete characterization of the nature of oscillations in system~\eqref{eq:nondimWel}, where each sub-region corresponds to a different combination in terms of the existence of  deep-(de)coupling phases. In particular, only in region $W$ one finds Welander oscillations with both a deep-coupling and a deep-decoupling phase. Note from the shape of region $W$ in Fig.~\ref{fig:BD}(a) that Welander oscillations require lower values of the density threshold $\eta$, with lower values of the virtual salinity flux $\mu$ for lower $\eta$. To illustrate further how the different types of oscillations arise, we present in panels~(b)--(d) of Fig.~\ref{fig:BD} the one-parameter bifurcation diagrams in $\mu$ for the three fixed values of $\eta$ that are indicated by the corresponding horizontal lines in panel~(a).

For $\eta = 0.05$ as in Fig.~\ref{fig:BD}(b), the periodic orbit $\Gamma$ appears and disappears for increasing $\mu$ at the Hopf bifurcations $\mathrm{H}$. In between these two points, the tangency curves $T^{\pm}$ are not encounterd. Therefore, $\Gamma$ remains entirely in the switching zone $S$ along the entire shown branch of periodic solutions, which is represented again by the maximum and minimum values of $\rho$ along $\Gamma$. Figure~\ref{fig:BD}(c) shows the bifurcation diagram in $\mu$ for $\eta = -0.17$, which is a value for which the tangency curves $T^{\pm}$ in panel~(a) are encountered, namely twice each. Immediately after the Hopf bifurcation that creates $\Gamma$ in Fig.~\ref{fig:BD}(c), the oscillation is confinced to the switching zone $S$. However, the tangency $T^{-}$ with $L^{-}$ takes place already for a slightly larger value of $\mu$. Subsequently, $\Gamma$ has a deep-decoupling phase, over the $\mu$-range up to the second encounter with the curve $T^{-}$ at $\mu \approx 0.3092$. The tangency $T^{+}$ with $L^{+}$ takes place only at $\mu \approx 0.1616$, and $\Gamma$ has a deep-coupling phase until $T^{+}$ is encountered again, which happens very close to the Hopf bifurcation $\mathrm{H}$ where $\Gamma$ disappears. Hence, there are Welander oscillations in the range $\mu \in [0.1616, 0.3092]$; the two tangencies shown in Fig.~\ref{fig:PPT} are exactly those at the two end points $T^{+}$ and $T^{-}$ of this range for $\eta = -0.17$.

In Fig.~\ref{fig:BD}(d) for the even lower value of $\eta = -0.3$, the order in which the curves $T^{\pm}$ are encountered remains unchanged; however, the two left-hand instances of $T^{-}$ and $T^{+}$ are now extremely close to the left-most Hopf bifurcation point $\mathrm{H}$ and, similarly, the right-hand points $T^{-}$ and $T^{+}$ are extremely close to the right-most point $\mathrm{H}$. This means that, for $\eta = -0.3$, Welander oscillations with definite deep-coupling and deep-decoupling phases exist, for all practical purposes, throughout the entire $\mu$-range of existence of the stable periodic orbit $\Gamma$. As the two-parameter bifurcation diagram in panel~(a) shows, this is also the case along horozontal lines for even lower $\eta$, as long as they are above the point $\mathrm{N}$.  Notice also that the size of $\Gamma$, which is represented by $\rho_{\min}$ and $\rho_{\max}$ in Fig.~\ref{fig:BD}(d), grows (or shrinks) extremely rapidly over tiny $\mu$-intervals near the two points $\mathrm{H}$ of Hopf bifurcation. This phenomenon is typical for slow-fast systems and known as a canard explosion \cite[]{slowfast_survey, kuehn2015multiple}.

\subsection{The four types of mixing oscillations}
\label{sec:phaseportraits}

\begin{figure*}[t!]
       \includegraphics{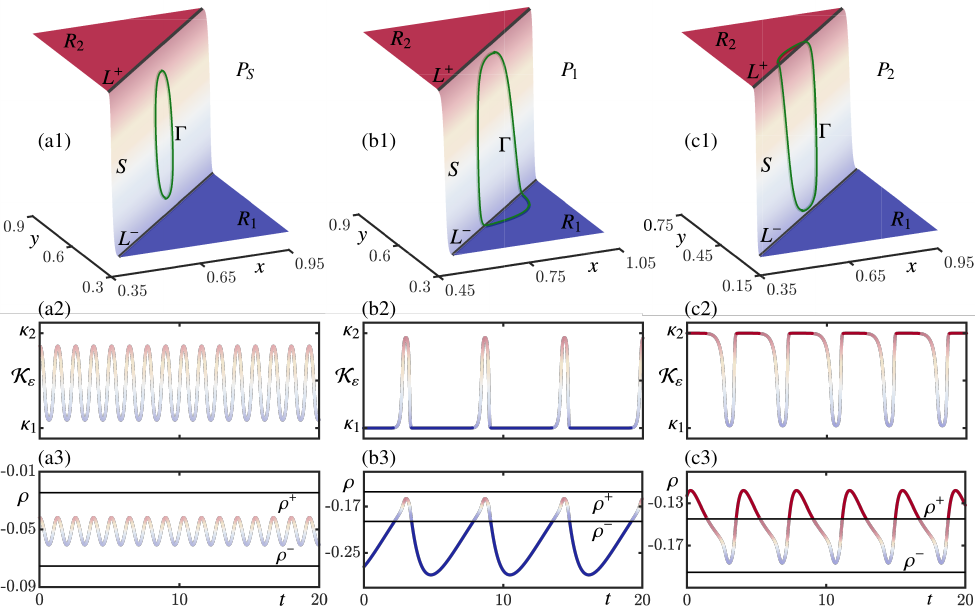}
	\caption{The periodic orbit $\Gamma$ of system~\eqref{eq:nondimWel} with $\varepsilon = 0.009$, when it is of type $P_S$ at $(\mu,\eta) = (0.32,-0.05)$ in panels~(a1)--(a3), of type $P_1$ at $(\mu,\eta) = (0.11, -0.17)$ in panels~(b1)--(b3), and of type $P_2$ at $(\mu,\eta) = (0.32, -0.17)$ in panels~(c1)--(c3), shown as in Figs.~\ref{fig:Welander} and~\ref{fig:PPT}.}
	\label{fig:PPP}
\end{figure*}

The example in Fig.~\ref{fig:Welander} of a Welander oscillation with $(\mu,\eta) = (0.14, -0.3)$ actually shows $\Gamma$ at the star in region $W$ of Fig.~\ref{fig:BD}(a). The oscillations in the parameter sub-regions $P_S$, $P_1$ and $P_2$ are shown in the same way in Fig.~\ref{fig:PPP}, as represented by $\Gamma$ at the stars in these corresponding sub-regions of the $(\mu,\eta)$-plane in Fig.~\ref{fig:BD}(a). These illustrations of $\Gamma$ on the graph of $\mathcal{K}_\varepsilon$, with the times series of $\mathcal{K}_\varepsilon$ and $\rho$, provide comprehensive insight into the nature of the different phases that define each type of oscillation.

Figure~\ref{fig:PPP}(a1) shows that the periodic orbit in sub-region $P_S$ is confined to the switching zone $S$; its time series of $\mathcal{K}_\varepsilon$ in panel~(a2) and of $\rho$ in panel~(a3) are almost sinosoidal and have a short period, which is due to $\Gamma$ evolving on the faster timescale associated with $S$. Indeed, the oscillation in $P_S$ lacks distinct deep-(de)coupling phases, meaning that there is no discernible deep water formation. Figure\ref{fig:PPP}(b1) shows $\Gamma$ in $P_1$, which is characterized by a deep-decoupled phase in zone $R_1$ and fast excursions into zone $S$, with the absence of a deep-coupling phase. The time series of $\mathcal{K}_\varepsilon$ in panel~(b2) shows a clear timescale separation: the deep-decoupled phase continues for quite some time, and is periodically interrupted by short bursts of the period of the oscillation in panel~(a2), with no distinct coupling phase. Panel~(b3) shows that, during the deep-decoupled phase, the density $\rho$ has a value below $\rho^{-}$ for long times. Hence, there is predominantly diffusive mixing in a stably stratified water column, with only brief intervals during which deep water production attempts to recover but fails. In sub-region $P_2$ one finds the converse: $\Gamma$ now has a deep-coupled phase in zone $R_2$ with fast excursions into zone $S$, but no deep-decoupled phase. This is illustrated in Figure\ref{fig:PPP}(c1) on the graph of $\mathcal{K}_\varepsilon$ with the associated time series in panel~(c2). The quite long deep-coupling phase with strong convective mixing is interrupted periodically by the water column attempting to re-adjust to a stable stratified state. However, this is unsuccessful and mixing returns to being convective as the periodic solution re-enters the deep-coupled zone $R_1$. Panel~(c3) illustrates that the density is generally sufficient to maintain an active deep water formation with brief lower-density episodes that are not causing a meaningful shutdown of convective mixing.

\section{Drifting virtual salinity flux}
\label{sec:drift}

We now study the influence of a steadily increasing freshwater influx in the North Atlantic, which corresponds to slowly decreasing virtual salinity 
flux, represented by $\mu$. The focus here is on an intermediate range of the density threshold $\eta$, where one encounters the sub-regions $P_1$, $W$ and $P_2$. Our analysis so far reveals that decreasing $\mu$ leads to progressively shorter deep-coupling phases and longer deep-decoupling phases. Hence, the deep water formation weakens and the AMOCs heat transport to the Northern Hemisphere diminishes in a negative convective feedback loop \cite[]{dijkstra2005nonlinear}. Furthermore, the deep water re-adjusts for short periods when the surface and deep-water boxes are decoupled, which  can lead to low-frequency variations in the AMOC. 

To model the deepwater formation process in this open system, we allow the virtual salinity $\mu$ to drift with time $t$ according to the linear ramp function
\begin{align}
\label{eq:mu_flux_drift}
\mu(t) = (1 - r t) \mu_{\rm start} + r t \mu_{\rm end} 
\end{align}
with rate (or slope of the ramp) $r$, which is chosen to be sufficiently small to ensure that system~\eqref{eq:nondimWel} with \eqref{eq:mu_flux_drift} adiabatically tracks the respective stable solutions between $\mu_{\rm start}$ and $\mu_{\rm end}$. In other words, $\mu(t)$ changes slowly enough so that only gradual transitions between different types of deep-(de)coupling oscillations are observed. This approach enables us to identify qualitatively different density oscillations in system~\eqref{eq:nondimWel} under a real-time gradual freshwater influx --- as could be observed in the modern climate due to an increased meltwater influx from the Greenland Ice Sheet into the Labrador and Nordic seas \cite[]{bamber2012recent, trusel2018nonlinear, sasgen2012timing}. 

\begin{figure*}[t!]
\centering
\includegraphics[scale=1]{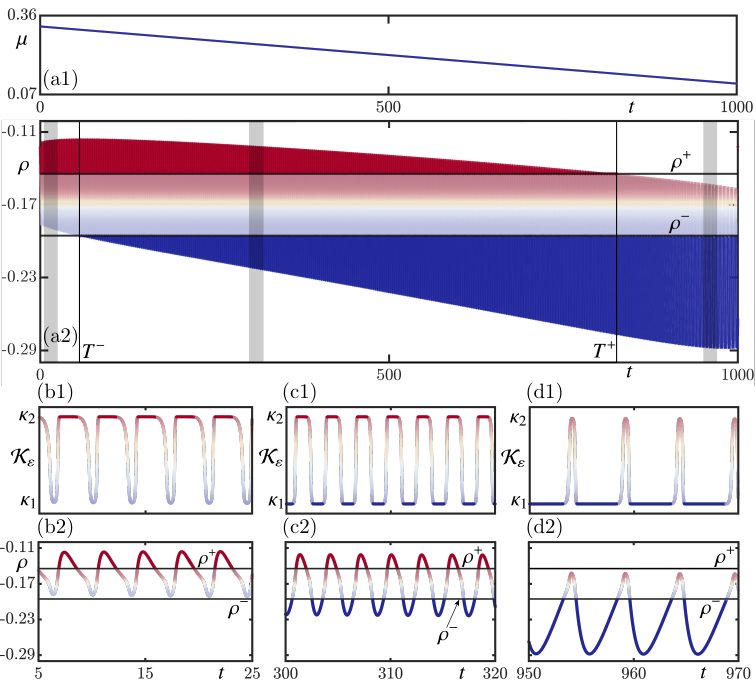}
\caption{Drifting dynamics of system \eqref{eq:nondimWel} with $\varepsilon = 0.009$ and for fixed $\eta = -0.17$ when the virtual salinity decreases linearly from $\mu_{\rm start} = 0.32$ to $\mu_{\rm end} = 0.11$ with rate of $r = 0.001$ in~\ref{eq:mu_flux_drift}. Panel~(a1) shows the ramping $\mu(t)$, and panel~(a2) the corresponding time series of  $\rho$ for $t \in [0, 1000]$. Also shown are the thresholds $T^{\pm}$ and three chosen time intervals of width $20$ (shaded in grey). The time series of $\mathcal{K}_\varepsilon$ in each of these time intervals are shown in panels~(b1), (c1) and~(d1), and those of $\rho$ in panels~(b2), (c2) and~(d2).}
\label{fig:nonauto_density}
\end{figure*}

Figure \ref{fig:nonauto_density} illustrates the transition as the virtual salinity flux $\mu(t)$ gradually decreases. It begins at the star in sub-region $P_2$ and ends at the star in sub-region $P_1$ from the bifurcation diagram in Figure \ref{fig:BD}(a); in particular, $\varepsilon = 0.009$ and $\eta = -0.17$ here. Specifically, we consider $\mu(t)$ as given by~\eqref{eq:mu_flux_drift}, with the start point $\mu_{\rm start} = 0.32$, as in Figure~\ref{fig:PPP}(c1)--(c3), and end point $\mu_{\rm end} = 0.11$, as Figure~\ref{fig:PPP}(b1)--(b3), and drift rate $r = 0.001$. This downward ramp is illustrated in Fig.~\ref{fig:nonauto_density}(a1) over the $1000$ time units it takes to complete the drift from start to finish. The corresponding time series of the density $\rho$ of system~\eqref{eq:nondimWel} is shown in panel~(a2), where the initial condition was a point on the periodic orbit $\Gamma$ in Figure~\ref{fig:PPP}(c1). Observe in Fig.~\ref{fig:nonauto_density}(a2) the overall decreasing trend of $\rho$, which is expected to accompany a freshening of the North Atlantic. This goes along with a shortening of the deep-coupling phase above $\rho^+$, and a considerable lengthening of the deep-coupling phase below $\rho^-$. Indeed, the former phase only exists until $\mu(t)$ drifts past the tangency $T^{+}$, and the latter only exists past $T^{-}$.

The drift rate $r$ is so small that the time series is effectively periodic over short time intervals; hence, the properties of the oscillation can be studied locally near a `frozen' value of $\mu$. We now examine the time series of $\mathcal{K}_\varepsilon$ and $\rho$ in the three windows, of 20 time units each, that are highlighted in Fig.~\ref{fig:nonauto_density}(a2). Panels (b1) and (b2) present these short time series in the first window with $t \in [5, 25]$, where we observe the characteristic properties of oscillations from sub-region $P_2$: long deep-decoupling phases interrupted by shorter dips in $\mathcal{K}_\varepsilon$ and $\rho$; compare with Figure~\ref{fig:PPP}(c2) and (c3). The second time window with $t \in [300, 320]$ is well in between $T^{-}$ and $T^{+}$, and the corresponding short time series in Fig.~\ref{fig:nonauto_density}(c1) and (c2) are indeed typical for Welander oscillations in sub-region $W$. They exhibit fast switching between prolonged deep-coupling and deep-decoupling phases of comparable length, which are here of about the same length here; compare with Figure~\ref{fig:Welander}(a2) and (a3). As $\mu(t)$ monotonically decreases, the length of the deep-decoupling phases increases, and the periodic convective flushes during the deep-coupling phases completely vanish when $T^{-}$ is passed. As is illustrated with Fig.~\ref{fig:nonauto_density}(d1) and (d2) for $t \in [950, 970]$, the short time series then clearly consists of long deep-coupling phases with periodic `blips' where the water column briefly looses its stratification, which is the hallmark of sub-region $P_1$; compare with Figure~\ref{fig:PPP}(b2) and (b3).

\section{Influence of the switching timescale parameter}
\label{sec:disapear_occ}

\begin{figure*}[t!]
\centering
\includegraphics{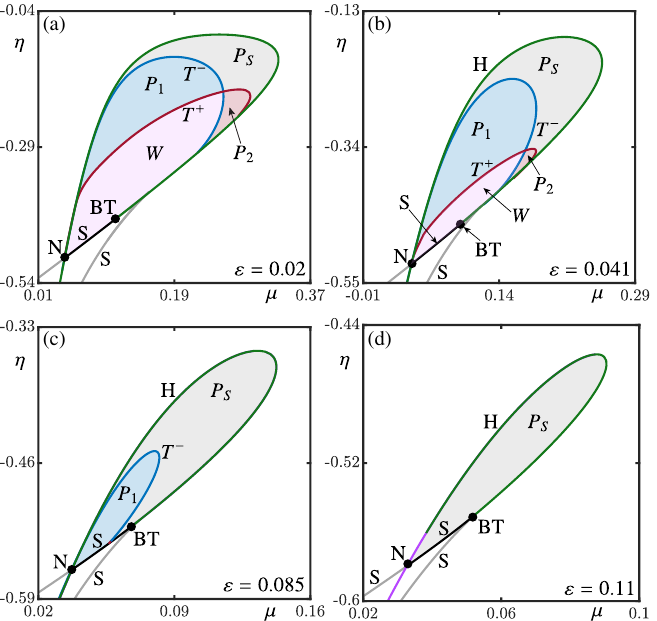}
\caption{Bifurcation diagram in the $(\mu,\eta)$-plane of system~\eqref{eq:nondimWel} with $\varepsilon = 0.02$ in panel~(a), $\varepsilon = 0.041$ in panel~(b), $\varepsilon = 0.085$ in panel~(c), and $\varepsilon = 0.11$ in panel~(d). Shown are the sub-regions $P_S$ (grey), $P_1$ (light blue), $P_2$ (light red) and $W$ (light purple), their bounding curves $T^{-}$ (blue), $T^{+}$ (red), $\mathrm{H}$ (green) and $\mathrm{S}$ (gray, and black when a SNIC), and the points $\mathrm{BT}$ and $\mathrm{N}$. Compare with Fig.~\ref{fig:BD}(a) and Fig.~\ref{fig:BDnonsmooth}(a).}
\label{fig:bd_disapearing_occ}
\end{figure*}

The parameter $\varepsilon$ directly influences the switching time between the convective and non-convective mixing phases, and it plays a crucial role in determining the configuration of deep-(de)coupling oscillations observed in system~\eqref{eq:nondimWel}. Since any oscillation completely disappears for $\varepsilon > 0.147$ \cite[]{bailie2024bifurcation}, we now address what this means for the fate of the sub-regions $P_S$, $P_1$, $P_2$, and $W$.

Figure~\ref{fig:bd_disapearing_occ} presents the bifurcation diagrams of system~\eqref{eq:nondimWel} in the $(\mu,\eta)$-plane for four increasingly larger values of $\varepsilon$. As before, the relevant parts of the curves $\mathrm{H}$ of Hopf bifurcation and $\mathrm{S}$ of saddle-node bifurcation bound the region of existence of the stable periodic solution $\Gamma$. Comparison of Fig.~\ref{fig:BD}(a) for $\varepsilon = 0.009$ with Fig.~\ref{fig:bd_disapearing_occ}(a) for $\varepsilon = 0.02$ and Fig.~\ref{fig:bd_disapearing_occ}(b) for $\varepsilon = 0.041$ reveals that these three bifurcation diagrams have the same qualitative features in terms of existence and locations of the sub-regions $P_S$, $P_1$, $P_2$ and $W$. In particular, they all feature Welander oscillations for lower values of $\eta$ above the point $\mathrm{N}$. We conclude that the situation in Fig.~\ref{fig:BD} is representative of any positive $\varepsilon$ over this initial range. In particular, as $\varepsilon$ is decreased to $0$ the Hopf bifucation curve $\mathrm{H}$, and the tangency curves $T^{-}$ and $T^{+}$ converge to the discontinuity indced bifurcations $\mathrm{BE}_1$, $\mathrm{FF}$ and $\mathrm{BE}_2$ in Figure~\ref{fig:BDnonsmooth}(a). In the process, the sub-regions $P_S$, $P_1$, and $P_2$ disappear, and only region $W$ of Welander oscillations remains in the limit $\varepsilon = 0$. 

While the bifurcation diagrams in Fig.~\ref{fig:BD}(a) and panels~(a) and (b) of Fig.~\ref{fig:bd_disapearing_occ} share their qualitative features, there are important quantitative differences. The existence region exhibiting oscillations decreases in size with increasing $\varepsilon$ (notice the different scales of the panels of Fig.~\ref{fig:bd_disapearing_occ}). At the same time, the relative sizes of the sub-regions $P_S$ and $P_1$ increase, while $P_2$ and $W$ decrease significantly; they are quite small in Fig.~\ref{fig:bd_disapearing_occ}(b) for $\varepsilon = 0.041$. As is shown in Fig.~\ref{fig:bd_disapearing_occ}(c),  when the switching timescale parameter is increased to $\varepsilon = 0.085$, the sub-regions $P_2$ and $W$ and their bounding curve $T^{+}$ have actually disappeared. This menas that the periodic orbit $\Gamma$ is always restricted to the switching zone $S$ and the deep-decoupling zone $R_1$. As a result, there are no longer any Welander oscillations, and deep water formation is suppressed. Note that the disappearance of Welander oscillations occurs substantially below $\varepsilon^{*} \approx 0.147$, where the region of oscillations itself vanishes. Figure~\ref{fig:bd_disapearing_occ}(d) presents the bifurcation diagram for $\varepsilon = 0.11$, where now also sub-region $P_1$ and its boundary curve $T^{-}$ have disappeared. Hence, the periodic solution always lies in the (now larger) switching zone $S$. 

Overall, we observe from Fig.~\ref{fig:bd_disapearing_occ} that the deep-coupling phase becomes increasingly limited with increasing $\varepsilon$. Already for `intermediate' values of $\varepsilon$, system \eqref{eq:nondimWel} is biased toward exhibiting deep-decoupling oscillations. Moreover, Welander oscillations, and deep-coupling phases more generally, only exist when the switching between zones $R_1$ and $R_2$ is sufficiently fast --- considerably faster than required for the existence of oscillations in the first place.

\section{Conclusions and outlook}
\label{sec:conclusions}

We have provided a detailed analysis of deep-decoupling oscillations in  the Welander model~\eqref{eq:nondimWel} with the smooth convective exchange function given by~\eqref{eq:convExchang}, which describes the transition between diffusive and convective mixing, the timescale of which is determined by the additional parameter $\varepsilon$. In the context of the deepwater  formation, this model exhibits self-sustained oscillations of temperature and salinity at the ocean surface. We conducted a bifurcation analysis that identifies sub-regions in the $(\mu, \eta)$-plane of virtual salinity influx and density threshold where one finds different types of deep-(de)coupling oscillations. The Welander oscillations with coexisting deep-decoupling and deep-coupling phases, which were shown to be the only type of oscillation for $\varepsilon = 0$, also exist for positive $\varepsilon$, but only up to a value considerable below the value $\varepsilon^{*}$ where oscillations disappear altogether. Indeed, this means that other types of oscillations exist for $\varepsilon >0$, but which exhibit only a deep-decoupling phase, only a deep-coupling phase, or neither of the two. We delimited their sub-regions of existence in the $(\mu, \eta)$-plane for any given $\varepsilon$ by computing curves of tangencies of the underlying stable periodic orbit with the boundary of the switching zone. Our bifurcation diagrams `extend' the known results for $\varepsilon = 0$ of instantaneous switching to the more general case where the switching is more gradual and described by a smooth convective exchange function with $\varepsilon > 0$. 

Overall, our results complete the analysis of the possible oscillatory behaviour exhibited by the Welander model. As we have shown, they may provide insights into the potential qualitative changes in the deepwater  formation under an influx of freshwater. Specifically, we found that large influxes result in a shortening of the convective phase, leading to prolonged adjustment during deep-decoupled phases. Conversely, with lower levels of freshening,  the deepwater  formation strengthened  due to longer phases of convective mixing and reduced weaker mixing. 

We distinguished the deep-decoupling, switching and deep-coupling zones by the lines where the curvature of the convective exchange function is maximal. This is convenient mathematical choice that works well for all $\varepsilon \geq 0$, but other choices for the boundaries of the switching zone could be made, such as using a percentage deviation from the minimum and maximum values $\kappa_1$ and $\kappa_2$. Moreover, one might consider a smooth switching function different from the hyperbolic tangent in~\eqref{eq:convExchang}. Whatever the choice, the technique of finding the curves of tangency with these zone boundaries can still be used to identify where the different types of oscillations exist. Moreover, if the alternative boundaries are close to the lines $L^{\pm}$ we considered here, then the respective sub-regions will only differ a little. Hence, what we presented is broadly representative for any reasonable choice of smooth switching function and switching thresholds. Since a form of convective adjustment is a common way to parameterize convective vertical mixing  in (particular relatively low  resolution) ocean models, our general approach can also be applied to the underlying convective adjustment scheme in more complex models. This would allow one to identify thresholds between deep-coupling and deep-decoupling regimes with realistic representations of the ocean and atmosphere dynamics. In particular, it will be interesting to determine the timescales of these two phases, in comparison with that of the (generally faster) switching between them; indeed, their ratio would give a `physical interpretation' to the switching timescale parameter $\varepsilon$ of the Welander model. In this way, it may also be possible to provide insights into the role of deep-decoupling oscillations in modelled D-O events  \cite[]{Vettoretti2022}. 

A next step in the vein of conceptual modelling would be to consider the competition between convective and salt-advective feedback mechanisms in the AMOC  \cite[]{weijer2019stability}. This can be achieved by adding additional boxes to the current schematic in Figure~\ref{fig:schematic}. Alternatively, one could incorporate an advective term with a time delay to represent the AMOCs adjustment period. However, the latter yields a model in the form of a delay differential equation with an infinite dimensional state space, the analysis of which comes with its own challenges \cite[]{kkp_dde_climate}.

\appendix
\section{Computing tangency curves}
\label{sec:computeT}

Equilibria, periodic orbits and their bifurcations can be found and followed in parameters with well-established numerical continuation techniques; see, for example, \cite[]{redbook2007} as an entry point to the literature. We use here the continuation software package {\sc AUTO} \cite[]{doedel2007auto} to find the stable periodic orbit $\Gamma$, as well as the bifurcation curves that bound the region where it exists; see \cite[]{bailie2024bifurcation} for more information on the exact nature of the bifurcations involved. 

The important new aspect here is the computation of the curves $T^{-}$ and $T^{+}$ by continuation, which we achieve as follows. The periodic orbit $\Gamma$ is continued from a point of Hopf bifurcation into the region $W$, where it has two intersetion points with each of the switching lines $L^{-}$ and $L^{+}$. The periodic solution $\Gamma$ is found as an orbit segment $\boldsymbol{u}$ that satisfies a boundary value problem (BVP); more specifically, $\boldsymbol{u}(t) = (\boldsymbol{u}_x(t), \boldsymbol{u}_y(t))$ satisfies the differential equation~\eqref{eq:nondimWel} (after rescaling with the period $T_\Gamma$ of $\Gamma$) subject to the periodicity conditions
\begin{equation}
\label{eq:po_cond}
\begin{aligned}
\begin{cases}
\boldsymbol{u}_x(0) =  \boldsymbol{u}_x(1),\\
\boldsymbol{u}_y(0) =  \boldsymbol{u}_y(1).
\end{cases}
\end{aligned}
\end{equation}
Note that the point $\boldsymbol{u}(0)$ could be any point along $\Gamma$ for conditions~\eqref{eq:po_cond} to be satisfied. This is why, to ensure this BVP is well defined and has a unique solution, one needs to impose a so-called phase condition \cite[Chapter~1 by E.~J.~Doedel]{redbook2007}.  {\sc AUTO} uses an integral phase condition \cite[]{doedel2007auto}, but we impose now the alternative condition
\begin{equation}
\label{eq:po_rho}
G(\boldsymbol{u}_y(0) - \boldsymbol{u}_x(0)) = 0,
\end{equation}
where the function $G$ is as defined in~\eqref{eq:max_curve_root}. Hence, the point $\boldsymbol{u}(0)$ is now a point on $L^{-}$ or $L^{+}$, and this also defines the solution $\boldsymbol{u}$ uniquely. 

We remark that an initial solution $\boldsymbol{u}$ satisfying~\eqref{eq:po_cond} and~\eqref{eq:po_rho} can be found in {\sc AUTO} by `rotating' the data representing $\Gamma$ until~\eqref{eq:po_rho} is satisfied. This orbit segement $\boldsymbol{u}$ can then be continued in a parameter, say, in $\mu$. When $\Gamma$ becomes tangent to the corresponding line $L^{-}$ or $L^{+}$, the continuation detects a fold point in the parameter (which is a double root of $G$). This fold can then be continued in an additional parameter, say, in $\eta$, yielding the respective curve $T^{-}$ and $T^{+}$ in the $(\mu,\eta)$-plane.

\bibliographystyle{elsarticle-num} 
\bibliography{BDK_AMOC_ARVIX.bib}

\end{document}